\documentclass{article}
\usepackage{graphicx} 
\usepackage{hyperref}

\title{A non-circular concept of number inspired by Gottlob Frege's definition}
\author{\Large{Prof. Marco Aurélio Spohn} \footnote{ \url{https://orcid.org/0000-0002-9265-9421}, ResearcherID: F-8993-2012, Scopus Author ID: 8435576500} \\ Federal University of Fronteira Sul \\ Chapecó, SC - Brazil}

\date{July 2024}

\begin{document}

\maketitle

\section{Introduction}

Gottlob Frege ingeniously presented a purely logical definition of the concept of number \cite{Frege1960,Frege2019}. However, one can claim that his definition is, in some way, circular, as it relies on the concept of one-to-one relation. When an object x is related to an object y, necessarily, y is related to x. In this sense, we assume the definition of the number one as a precondition.

The definition presented in this paper humbly borrows the description style from Frege's work. However, we emphasize that this work falls far short of a description as formal as that of Frege's original work. The original inspiration and all thanks are due to the work of such a brilliant philosopher and mathematician.

\section{A non-circular definition of number}

The concept of number only makes sense when it presents the property of projection/reflection or binding. When we consider a number as an abstraction of objects, whatever they may be, saying that a number that belongs to the concept F is the same as that which belongs to the concept G means there is a projection/reflection, or binding, between the objects in F and the objects in G.

We present a definition based on both equivalent approaches. First, we introduce the definition based on the relations of projection and reflection; then, we present the definition based on the relation of binding.

\subsection{Definition based on the projection/reflection relations}

At first glance, when thinking about the concept of numbers, the idea of comparing some accurate or abstract representation of objects with a reference concept arises. Then, the first impression of this process comes with the relation of projection; that is, objects falling under a given concept are projected and, consequently, reflected falling under another concept.

\paragraph{\textbf{Projection}} can be defined as the relation between objects \textbf{x} and \textbf{y} when only object \textbf{x} projects onto object \textbf{y}, excluding any object other than \textbf{x}.

\paragraph{\textbf{Reflection}} can be defined as the relation between objects \textbf{y} and \textbf{x} when only object \textbf{y} reflects the projection of \textbf{x}, excluding any object other than \textbf{y}.

From the projection and reflection relations, the following can be deemed as well understood:
\begin{itemize}
\item If \textbf{a} stands in the projection relation to \textbf{b}, and \textbf{a} stands in the projection relation to \textbf{c}, then, generally, whatever \textbf{a}, \textbf{b}, and \textbf{c} are, \textbf{b} and \textbf{c} are the same.
\item If \textbf{b} stands in the reflection relation to \textbf{a}, and \textbf{d} stands in the reflection relation to \textbf{a}, then, generally, whatever \textbf{a}, \textbf{b}, and \textbf{d} are, \textbf{b} and \textbf{d} are the same.
\end{itemize}

There is a relation $\phi$ that projects all objects falling under concept F onto the objects falling under concept G unambiguously; that is, every object falling under concept F has exclusivity in the projection onto some object falling under concept G and, likewise, every object falling under concept G has exclusivity in the projection onto some object falling under concept F. The exclusivity results from the mutual correspondence between the projection and reflection relations.

It is worth mentioning that at no moment is it implied that there is a strict \textit{one-to-one relation}\footnote{In a one-to-one relation/correspondence (bijection) between objects from sets A and B, if an object \textbf{a} from set A is mapped to an object \textbf{b} from set B, then from the perspective of set B, the object \textbf{b} is necessarily mapped back to the object \textbf{a} from set A. This is the essence of a bijective relationship — a two-way, reciprocal mapping.} between objects falling under concept F and objects falling under concept G: projections/reflections to/from objects falling under concept F and falling under concept G are independent. 

\begin{figure}
    \centering
    \includegraphics[width=0.85\columnwidth]{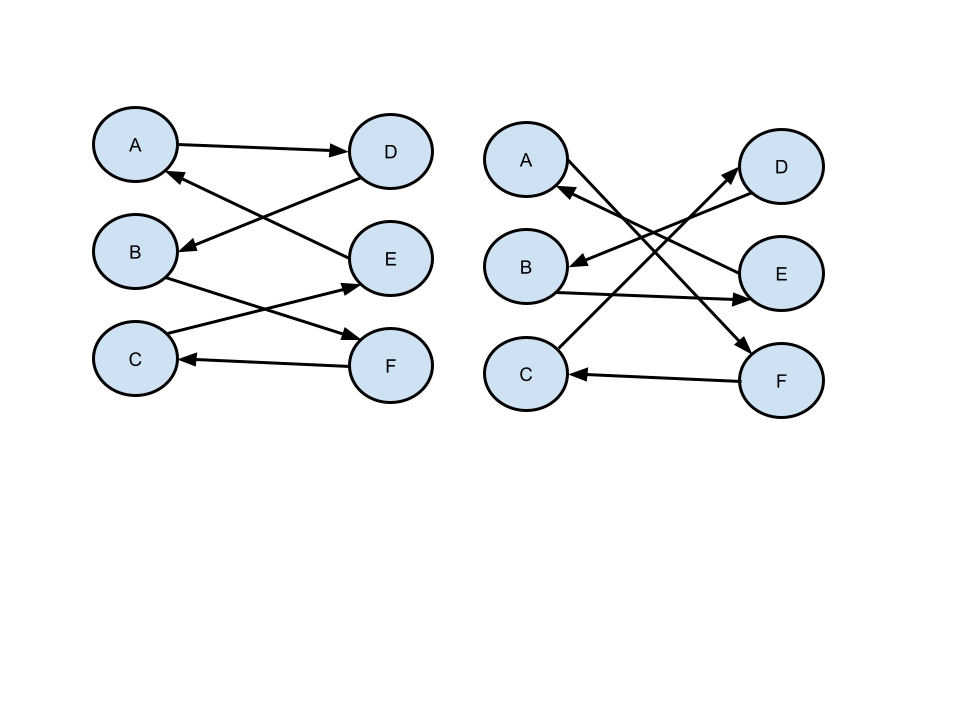}
    \caption{Projection or binding: diversity arising from multiple valid combinations.}
    \label{fig:relation}
\end{figure}

\paragraph{The number} that belongs to concept F is the extension of the concept `projection/reflection to the concept F.'

\paragraph{The expression}  ``\textbf{n} is a number'' is synonymous with the expression  ``There is a concept such that \textbf{n} is the number that belongs to it.''.

Any possibility of projection and reflection resulting from the relation $\phi$ is valid; any possible combination reflects (refers to)  the same number (see Figure~\ref{fig:relation}).

\paragraph{Regarding the definition of number one,} it is obtained directly from our general definition. However, one could also define it by saying that number one is an object that projects onto itself and, consequently, projects itself.

\paragraph{As for the number zero,} it can be defined as an object that projects onto no object and, consequently, has no reflection.

\subsection{Definition based on the binding relation}

From the perspective of binding, there are the following relations:

\paragraph{The binding relation, \textit{stringTo},} between objects \textbf{x} and \textbf{y} occurs when object \textbf{x} has exclusive binding to object \textbf{y}.
 
\paragraph{The binding relation, \textit{stringFrom},} between objects \textbf{y} and \textbf{x} occurs when object \textbf{y} is in the binding from object \textbf{x}.

From the binding relations, \textit{stringTo} and \textit{stringFrom}, the following is well understood:

\begin{itemize}
\item If object \textbf{a} stands in the \textit{stringTo} relation to object \textbf{b}, and object \textbf{a} stands in the \textit{stringTo} relation to object \textbf{c}, then, generally, whatever \textbf{a}, \textbf{b}, and \textbf{c} are, \textbf{b} and \textbf{c} are the same.
\item If object \textbf{b} stands in the \textit{stringFrom} relation to object \textbf{a}, and object \textbf{d} stands in the \textit{stringFrom} relation to object \textbf{a}, then, generally, whatever \textbf{a}, \textbf{b}, and \textbf{d} are, \textbf{b} and \textbf{d} are the same.
\end{itemize}

There is a relation $\phi$ that binds all objects falling under concept F with the objects falling under concept G unambiguously; that is, every object falling under concept F has exclusivity in the binding to some object falling under concept G and, likewise, every object falling under concept G has exclusivity in the binding to some object falling under concept F. The exclusivity is given by the mutual correspondence between the \textit{stringTo} and \textit{stringFrom} relations.

Once again, it is worth mentioning that at no moment is it implied that there is a strict one-to-one relation between objects falling under concept F and objects falling under concept G: bindings to/from objects falling under concept F and falling under concept G are independent. 

\paragraph{The number} that belongs to concept F is the extension of the concept `binding to the concept F.'

\paragraph{The expression}  ``\textbf{n} is a number'' is synonymous with the expression ``There is a concept such that \textbf{n} is the number that belongs to it.''.

Any possibility of binding resulting from the relation $\phi$ is valid; any possible binding configuration relates to the same number.

\paragraph{Regarding the definition of number one,} it is obtained directly from our general definition. However, one could also define it by saying that the number one is an object that binds to itself.

\paragraph{As for the number zero,} it can be defined as an object that binds to no object.

\section{Conclusion}

Frege's concept of numbers, presented in 1884, represents a phenomenal contribution to the philosophy of mathematics. Although it is argued that his definition is somehow circular, its unique and utterly logical presentation makes it the primary reference to the definition of number.

This work presents two approaches to defining numbers in a non-circular way. We defend that numbers derive essentially from the relationship of projection and reflection, or binding, inherent to their existence. Considering these relationships, numbers can be defined without incurring circular definitions. The properties of relationships guarantee that there is no need to rely on the particular case of one-to-one associations. A more formal and detailed description of the foundations presented in this work is planned in future work.

\bibliographystyle{abbrv}
\bibliography{references}

\end{document}